\catcode`\@=11

\magnification=1200
\baselineskip=14pt

\pretolerance=500    \tolerance=1000 \brokenpenalty=5000

\catcode`\;=\active
\def;{\relax\ifhmode\ifdim\lastskip>\z@
\unskip\fi\kern.2em\fi\string;}

\overfullrule=0mm

\catcode`\!=\active
\def!{\relax\ifhmode\ifdim\lastskip>\z@
\unskip\fi\kern.2em\fi\string!}

\catcode`\?=\active
\def?{\relax\ifhmode\ifdim\lastskip>\z@
\unskip\fi\kern.2em\fi\string?}

\frenchspacing

\newif\ifpagetitre            \pagetitretrue
\newtoks\hautpagetitre        \hautpagetitre={ }
\newtoks\baspagetitre         \baspagetitre={1}

\newtoks\auteurcourant        \auteurcourant={M. Laurent   }
\newtoks\titrecourant
\titrecourant={ Exponents of Diophantine approximation in dimension two }

\newtoks\hautpagegauche       \newtoks\hautpagedroite
\hautpagegauche={\hfill\sevenrm\the\auteurcourant\hfill}
\hautpagedroite={\hfill\sevenrm\the\titrecourant\hfill}

\newtoks\baspagegauche       \baspagegauche={\hfill\rm\folio\hfill}

\newtoks\baspagedroite       \baspagedroite={\hfill\rm\folio\hfill}

\headline={
\ifpagetitre\the\hautpagetitre
\global\pagetitrefalse
\else\ifodd\pageno\the\hautpagedroite
\else\the\hautpagegauche\fi\fi}

\footline={\ifpagetitre\the\baspagetitre
\global\pagetitrefalse
\else\ifodd\pageno\the\baspagedroite
\else\the\baspagegauche\fi\fi}

\def\date{\ {\the\day}\
\ifcase\month\or Janvier\or F\'evrier\or Mars\or Avril
\or Mai \or Juin\or Juillet\or Ao\^ut\or Septembre
\or Octobre\or Novembre\or D\'ecembre\fi\
{\the\year}}

\def\up#1{\raise 1ex\hbox{\sevenrm#1}}

\def\cqfd{\unskip\kern 6pt\penalty 500
\raise -2pt\hbox{\vrule\vbox to 10pt{\hrule width 4pt
\vfill\hrule}\vrule}\par\medskip}

\def\section#1{\vskip 7mm plus 20mm minus 1.5mm\penalty-50
\vskip 0mm plus -20mm minus 1.5mm\penalty-50
{\bf\noindent#1}\nobreak\smallskip}

\def\subsection#1{\medskip{\bf#1}\nobreak\smallskip}

\def\displaylinesno #1{\dspl@y\halign{
\hbox to\displaywidth{$\@lign\hfil\displaystyle##\hfil$}&
\llap{$##$}\crcr#1\crcr}}

\def\ldisplaylinesno #1{\dspl@y\halign{
\hbox to\displaywidth{$\@lign\hfil\displaystyle##\hfil$}&
\kern-\displaywidth\rlap{$##$}
\tabskip\displaywidth\crcr#1\crcr}}

\def\hfl#1#2{\smash{\mathop{\hbox to 12 mm{\rightarrowfill}}
\limits^{\scriptstyle#1}_{\scriptstyle#2}}}

\catcode`\@=12

\def\bP{{\bf P}}
\def\bQ{{\bf Q}}
\def\bR{{\bf R}}
\def\bZ{{\bf Z}}

\def\uA{{\underline{A}}}
\def\uB{{\underline{B}}}
\def\uD{{\underline{D}}}

\def\uP{{\underline{P}}}

\def\uQ{{\underline{Q}}}

\def\uX{{\underline{X}}}

\def\uDelta{{ \underline{\Delta}}}

\def\bZ{{\bf Z}}

\def\proof{\bigskip\noindent{\it Proof.}\ }

\def\and{\quad\hbox{and}\quad}

\def\om{{\omega}}
\def\omc{{\hat{\omega}}}


\def\M{\mathop{\rm M\kern 1pt}\nolimits}
\def\h{\mathop{\rm h\kern 1pt}\nolimits}

\def\romain#1{\uppercase\expandafter{\romannumeral #1}}

\def\card{\mathop{\rm Card\kern 1.3 pt}\nolimits}
\def\deg{\mathop{\rm deg\kern 1pt}\nolimits}
\def\det{\mathop{\rm det\kern 1pt}\nolimits}

\def\h{\mathop{\rm h\kern 1pt}\nolimits} \long\def\forget#1\endforget{}

\def\og{\leavevmode\raise.3ex\hbox{$\scriptscriptstyle 
\langle\!\langle\,$}}
\def\fg{\leavevmode\raise.3ex\hbox{$\scriptscriptstyle
\!\rangle\!\rangle\,\,$}}

\def\Bak{1}
\def\BoVa{2}
\def\BuLaA{3}
\def\BuLaB{4}
\def\BuLaC{5}
\def\DaSc{6}
\def\Fis{7}
\def\JarA{8}
\def\JarB{9}
\def\JarC{10}
\def\JarD{11}
\def\JarE{12}
\def\JarF{13}
\def\Khi{14}
\def\Lag{15}
\def\RoyA{16}
\def\RoyB{17}
\def\RoyC{18}
\def\Ry{19}
\def\Sch{20}

\centerline{}

\vskip 4mm

\centerline{
\bf  Exponents of  Diophantine approximation in dimension two}

\vskip 8mm
\centerline{ by Michel L{\sevenrm AURENT}\footnote{}{\rm
2000 {\it Mathematics Subject Classification : }  11J13, 11J70.} }

\vskip 9mm

\noindent {\bf Abstract --} Let $\Theta=(\alpha,\beta)$ be a point in $\bR^2$, with $1,\alpha,\beta$
linearly independent over $\bQ$.
We attach to $\Theta$ a quadruple $\Omega(\Theta)$ of  exponents which measure  the quality
of approximation to $\Theta$ both by rational points and by rational lines.
The two ``uniform'' components of $\Omega(\Theta)$ are related by an equation, due to Jarn\'\i k, 
and the four exponents satisfy  two inequalities which refine Khintchine's transference principle.
Conversely, we show that for any quadruple $\Omega$ fulfilling  these necessary  conditions, there exists
 a point $\Theta\in \bR^2$ for which $\Omega(\Theta) =\Omega$. 

\vskip15mm

\section
{1. Introduction and results.}
Let $\alpha$ and $\beta$ be real numbers. 
  We  first introduce     four  exponents which quantify various notions of rational approximation to the point $(\alpha,\beta)$
  in the plane $\bR^2$. 
  
Define     $\om(\alpha,\beta)$  as    the  supremum (possibly infinite)  
of  all real numbers $\om$ such that  {\it there exist infinitely many integers} $H$ for which the 
 inequalities
$$
| x\alpha +y \beta +  z | \le H^{-\om} \and \max \{ | x|,|y|,| z |\} \le H
$$
admit a non-zero integer solution $(x,y,z)$.  Following the  general notations  of [\BuLaB], we define   moreover   $\omc(\alpha,\beta)$  
as  the supremum of  all real numbers $\om$ such that for {\it any sufficiently large integer} $H$, 
the above system of inequations has a non-zero integer solution.
Considering as well the simultaneous rational approximation to $\alpha$ and $\beta$, 
we  define  similarly two further exponents $\om\left(\matrix{\alpha\cr\beta\cr}\right)$ and 
$\omc\left(\matrix{\alpha\cr\beta\cr}\right)$ by repeating word for word the previous sentences,   and  replacing  the above  inequalities by 
$$
\max \left\{ | z\alpha - x  |,|z\beta -y| \right\} \le H^{-\om} \and \max \{ | x|,|y|,| z |\} \le H.
$$

 The  exponents $\om(\alpha,\beta)$ and $\om\left(\matrix{\alpha\cr\beta\cr}\right)$
 are those which occur  most frequently in Diophantine approximation. 
 Substituting    $\max \{ | x|,|y|,| z |\}$ for $H$ 
in the preceding  inequations,  we observe  that  these two exponents   measure  respectively the sharpness 
of the approximation to the point $(\alpha,\beta)$ by  rational lines and by  rational points, in terms of  their height. 
 The   corresponding uniform exponents 
$\omc(\alpha,\beta)\ge 2$ and $\omc\left(\matrix{\alpha\cr\beta\cr}\right)\ge 1/2$ were first  introduced  by Jarn\'\i k.
They quantify the possible improvements to Dirichlet  box principle when applied to the two  systems of linear 
inequalities. 

Set $\Theta =(\alpha,\beta)$ and denote by
${}^t\Theta=\left(\matrix{\alpha\cr\beta\cr}\right)$ the transposed vector. For brevity, we shall often write
$$
\om(\Theta)=\om(\alpha,\beta),\quad \om({}^t\Theta)= \om\left(\matrix{\alpha\cr\beta\cr}\right), \quad 
\omc(\Theta) = \omc(\alpha,\beta), \quad \omc({}^t\Theta)=\omc\left(\matrix{\alpha\cr\beta\cr}\right).
$$
The goal of our article is to describe the {\it spectrum } of these four exponents, that is the set of values
 taken by the quadruples 
 $$
 \Omega(\Theta) = (\om(\Theta),\om({}^t\Theta), \omc(\Theta),\omc({}^t\Theta)),
 $$
 when $\Theta=(\alpha,\beta)$ ranges over $\bR^2$, with $1,\alpha,\beta$ linearly independent over $\bQ$.
We  have conventionally excluded from the spectrum the   points  with  $1,\alpha,\beta$
  linearly dependent over $\bQ$, for which the four exponents behave as for real numbers. 
In this latter case, observe that
 the   exponent $\omc\left(\matrix{\alpha\cr 0\cr}\right)$ of uniform rational approximation  to $\alpha$ 
is equal to $1$,  whenever  $\alpha$  is irrational (Satz 1 of Khintchine's seminal paper [\Khi]). 
Thus, if  the numbers $1,\alpha,\beta$ are linearly dependent over $\bQ$ and at least one of  the numbers $\alpha$ or 
$\beta$ is irrational, the quadruple $\Omega(\Theta)$ has the form 
$$
\Omega(\Theta) = (+\infty ,v, +\infty , 1)
$$
with $v\ge 1$, and any value  $v$ in the interval $[1,+\infty]$ may be reached for some  point $\Theta$. 
When both  $\alpha$ and $\beta$ are rational, we obviously have
 $$
 \Omega(\Theta) = (+\infty, +\infty, +\infty, +\infty).
 $$
  From now, we shall assume that the numbers $1,\alpha,\beta$ are linearly independent over $\bQ$.

 Jarn\'\i k has studied the relations between  the exponents $\om$ and  $\omc$  in a series of papers [\JarD--\JarF]
 dealing with any system of real linear forms. 
 We refer to [\BuLaB --\BuLaC] for a detailled  survey of his results  on this topic. In dimension two,  he proved  [\JarD]  the formula
 $$
\omc\left(\matrix{\alpha\cr\beta\cr}\right) ={\omc(\alpha,\beta)\over \omc(\alpha,\beta)-1}.
$$
The exponents $\om(\alpha,\beta) $ and $\om\left(\matrix{\alpha\cr\beta\cr}\right)$ are
related by  Khintchine's transference inequalities 
$$
{\om(\alpha,\beta)\over \om(\alpha,\beta)+2} \le \om\left(\matrix{\alpha\cr\beta\cr}\right) \le {\om(\alpha,\beta) -1\over 2}.
$$
See for instance Satz VI of [\Khi]. Our theorem refines this latter estimate.

\proclaim
Theorem.
For any row vector  $\Theta = (\alpha, \beta) $ with $1,\alpha , \beta$ linearly independent over $\bQ$,  the four exponents
$$
v= \om(\Theta),\quad  v'= \om({}^t\Theta),
\quad  w = \omc(\Theta), \quad w'= \omc({}^t\Theta), 
$$
satisfy  the relations
$$
2\le w\le +\infty ,\quad  w' ={w-1\over w} , \quad {v(w-1)\over v+w} \le v' \le {v-w+1\over w}.
$$
When  $w < v=+\infty$ we have to understand these relations as   $w-1 \le v' \le +\infty$,  and when $w=+\infty$, the set of constraints
should be interpreted as $v=v'=+\infty$ and $w'=1$.
Conversely,  for each quadruple  $(v,v',w,w')$ in  $(\bR_{>0}\cup\{+\infty\})^4$ 
satisfying the previous conditions, 
there exists a row vector $\Theta=(\alpha,\beta)$ of real numbers with  $1,\alpha , \beta$ linearly independent over $\bQ$,
such that
$$
\Omega(\Theta) = (v,v',w,w').
$$

Notice that the  estimate
$$
{v(w-1)\over v+w} \le v' \le {v-w+1\over w}
$$
refines Khintchine's  inequalities since $w\ge 2$. 

Few  explicit computations
of quadruples $\Omega(\Theta)$ 
have actually been achieved. It follows from   Roy's works  [\RoyA--\RoyB]  that 
$$
\omc(\alpha, \alpha^2) ={ 3+ \sqrt{5}\over 2}  , \quad  \omc\left(\matrix{\alpha\cr\alpha^2\cr}\right) = {  \sqrt{5}-1\over 2}, 
$$
when $\alpha$ is a so-called Fibonacci continued fraction.   Next, Bugeaud \& Laurent [\BuLaA] have explicitly determined the quadruple 
$\Omega(( \alpha,\alpha^2))$ for any  
  sturmian continued fraction $\alpha$. Further (very partial) informations
  on  quadruples   of  the form $\Omega( (\alpha,\alpha^2))$, where $\alpha$ is a
   real transcendental number, may  also be derived from   [\BuLaC , \Fis, \RoyC].

\medskip

Jarn\'\i k [\JarE--\JarF] has improved  the obvious lower bounds $\om(\Theta)\ge \omc(\Theta) $ and  $\om({}^t\Theta)\ge \omc({}^t\Theta)$.
We deduce his results from our theorem and we show that they are optimal.

\proclaim
Corollary 1. 
For any row vector  $\Theta =(\alpha, \beta)$ with $1,\alpha,\beta$ linearly independent over $\bQ$, 
 the lower   bounds 
$$
 \omc(\Theta) \ge 2 \and \om(\Theta) \ge \omc(\Theta)(\omc(\Theta)-1)
$$
hold. Conversely, for any  $v\in \bR_{>0}\cup\{+\infty\}$ and any $w\in \bR_{>0}\cup\{+\infty\}$    satisfying
$$
2 \le w \le +\infty     \and  w(w-1) \le v \le +\infty,
$$
  there exists a row vector $ \Theta =(\alpha, \beta)$ with $1,\alpha,\beta$ linearly independent over $\bQ$, such that
$$
\om(\Theta)=v \and \omc(\Theta) =w.
$$

\proclaim
Corollary 2. 
For any  column vector $\Theta =\left(\matrix{\alpha\cr\beta\cr}\right)$
 with $1,\alpha,\beta$ linearly independent over $\bQ$, we have 
 $$
{1\over 2} \le \omc(\Theta) \le 1 \and \om(\Theta) \ge {\omc(\Theta)^2\over 1 - \omc(\Theta)}.
$$
Conversely,  for any   $w' \in \bR_{>0}$ and any $v'\in\bR_{>0}\cup\{+\infty\}$  satisfying 
$$
{1\over 2 }\le w'\le 1 \and   { w'^2\over 1-w'} \le v' \le + \infty,
$$
 there exists a column vector $\Theta =\left(\matrix{\alpha\cr\beta\cr}\right)$
 with $1,\alpha,\beta$ linearly independent over $\bQ$, such that
$$
\om(\Theta)=v' \and \omc(\Theta) =w'.
$$

The existence of a column (resp. row) vector $\Theta$ for which $\omc(\Theta)$ takes an arbitrary
value   in the interval $[1/2,1]$ (resp. $[2,+\infty]$)  follows from [\JarF].  Jarnik's   approach, which is  based on some
explicit construction of continued fractions, differs from ours. 

In order to derive  both corollaries from the theorem, observe that, for given positive real numbers
$v$ and $w$, the interval
$$
{v(w-1)\over v+w} \le v' \le {v-w+1\over w}
$$
occurring in our theorem, is non-empty exactly when $v\ge w(w-1)$. 
For the minimal value $v=w(w-1)$, it reduces to the
point  
$$
{(w-1)^2\over w}= {w'^2\over 1-w'}.
$$ 
 Corollaries 1 and 2 immediately follow,  noting that the extremal  values  $v(w-1)/(v+w)$ and $(v-w+1)/w$
 are increasing functions of $v$, when $v\ge w(w-1)$. 
 
 \medskip 
 
 The proof of our theorem splits into two parts. We first establish  the two transference 
 inequalities by means of simple geometrical constructions involving the best rational
  approximations (``minimal points'' in the terminology of Davenport \& Schmidt [\DaSc]) to the point $\Theta$. 
  The determination of a point $\Theta$ with prescribed $\Omega(\Theta)$, 
    needs more elaborate arguments. We simultaneously construct  a sequence
  of rational lines $\Delta_{n,k}$ and a Cauchy sequence of rational  points $P_{n,k}$,  which approximate
  the limit  $\Theta= \lim P_{n,k}$ in a controlled way. The geometrical configuration of these two sequences
  of lines and of points (colinear points and concurrent lines) reflects  duality  relations
  between two  sequences of best approximations by lines and by points to a given point $\Theta\in \bR^2$.
  
  To conclude this introduction, let us address the problem of extending   the theorem in higher dimensions.
  Then $\Theta$ should stand for any real linear proper subvariety of a projective space $\bP^m(\bR)$,  to which 
  we  can attach various (usual and uniform) exponents of approximation by rational linear subvarieties
  of fixed dimension $\mu, 0\le \mu \le m-1$, as in [\BuLaC, \Sch]. 
   We refer to Section 4 of  [\BuLaC] for  precise definitions and  ask for a description of the  
   spectrum determined by the vector $\Omega(\Theta)$
  of  these exponents,  when $\Theta$ ranges over the set of all real linear subvarieties of $\bP^m(\bR)$ with given dimension. 
 As a next step after the present situation dealing with a point $\Theta$ in $\bP^2(\bR)$,
  it should be interesting to investigate the case of a point   in $\bP^3(\bR)$
   which gives rise to six exponents. 
   
   \bigskip
   \noindent
   {\bf Acknowledgements.} I would like to thank Damien Roy for his deep analysis of the paper which led us to improve 
   an earlier version of Lemmas 1 and 2.

\section
{2. Transference inequalities.}

We prove  in this section   the transference inequalities
$$
{v(w-1)\over v+w} \le v' \le {v-w+1\over w}
$$
  for any point $\Theta=(\alpha, \beta)$ with $1,\alpha,\beta$ linearly independent over $\bQ$.
Two specific sequences of best approximations will serve our purpose.
 We refer to [\BuLaB,  \Lag] for further informations on the notion of  best approximation.
 
For any  triple $\uX =(x,y,z)$ of real numbers, set
$$
L(\uX) = | x\alpha + y\beta +z|, \quad M(\uX) = \max (|z\alpha -x|,|z\beta -y|) \and \| \uX \| = \max\{| x|,| y|,| z|\}.
$$
We call  {\it sequences of best approximations} relative to the norm $\| \,\, \|$,   respectively associated to 
the semi-norms $L$ and $M$,   
 two sequences of integer triples
$$
\uDelta_n =(r_n,s_n,t_n)  \and \uP_n =(a_n,b_n, c_n) , \quad n\ge 1,
$$
satisfying   the following properties. Put
$$
h_n =\| \uDelta_n\| , \quad q_n = \| \uP_n\|, \quad L_n = L(\uDelta_n) , \quad M_n = M(\uP_n).
$$
The sequences of norms 
$$
1 <   h_1 < h_2 < \dots  \and \quad 1 <   q_1 < q_2 < \dots , 
$$
increase,  while  the sequences of values 
$$
1 > L_1> L_2 > \dots   \and  \quad 1 > M_1> M_2 > \dots . 
$$
decrease and tend to $0$. 
For any $n\ge 1$ and any non-zero  integer triple $\uDelta$ (resp. $\uP$) with norm $\|\uDelta\| < h_{n+1}$
(resp. $\| \uP\| < q_{n+1}$), we have the lower bounds
$$
L(\uDelta) \ge L_n \and M(\uP) \ge M_n.
$$

Define now  positive exponents $v_n, v'_n,w_n, w'_n$ by the equations
$$
L_n = h_n^{-v_n}= h_{n+1}^{-w_n} \and M_n = q_n^{-v'_n} =q_{n+1}^{-w'_n} , \quad (n\ge 1).
$$
Our  interest in these two sequences of best approximations  rests on 
the formulas
$$
\om(\Theta)= \limsup_{n\rightarrow +\infty} v_n , \quad \om({}^t\Theta)= \limsup_{n\rightarrow +\infty} v'_n , 
\quad \omc(\Theta) = \liminf_{n\rightarrow +\infty} w_n , 
 \quad \omc({}^t\Theta) = \liminf_{n\rightarrow +\infty} w'_n ,
$$
which  enable us to compute  $\Omega(\Theta)$  thanks to  the sequences of test  points
$\uDelta_n$ and  $\uP_n$, as it is easily seen from the above  properties.  

A   geometrical point of view may be enlightening. Denote by $\Delta_n$ the line in $\bR^2$
with equation $r_nx+s_ny +t_n=0$,  and by $P_n$ the rational point with coordinates $P_n = (a_n/c_n,b_n/c_n)$.
In the sequel we shall follow  these  conventions of notations.  An underlined symbol will always stand for some non-zero real triple.
The same  symbol without underlining  will indicate either the associated  line (as for $\Delta_n$), 
or the point obtained by deshomogenization with respect to the third coordinate (as for $P_n$).  
The alternative will be clear from the context.

Observe now that two consecutive best approximations $\uDelta_n$ and $\uDelta_{n+1}$ are not proportional. Therefore the
vector product
$$
\uQ_n = \uDelta_n\wedge \uDelta_{n+1} = \left( \left\vert \matrix{s_n & s_{n+1}\cr t_n & t_{n+1}\cr}\right\vert , 
- \left\vert \matrix{r_n & r_{n+1}\cr t_n   & t_{n+1}\cr}\right\vert , 
\left\vert \matrix{r_n & r_{n+1}\cr s_n & s_{n+1}\cr}\right\vert \right), 
$$
is a non-zero triple, so that   $\Delta_n$ cuts  $\Delta_{n+1}$ at the point $Q_n$. 
Since both  lines $\Delta_n$  and $\Delta_{n+1}$ are close to $\Theta$, their intersection $Q_n$ should  also be close to $\Theta$.
More precisely, write 
$$
\eqalign{
\left\vert \matrix{r_n & r_{n+1}\cr s_n & s_{n+1}\cr}\right\vert \alpha -  \left\vert \matrix{s_n & s_{n+1}\cr t_n & t_{n+1}\cr}\right\vert 
=  & s_{n+1}(r_n\alpha +s_n \beta +t_n) - s_n(r_{n+1}\alpha +s_{n+1}\beta +t_{n+1}) ,
\cr
\left\vert \matrix{r_n & r_{n+1}\cr s_n & s_{n+1}\cr}\right\vert \beta  +   \left\vert \matrix{r_n & r_{n+1}\cr t_n & t_{n+1}\cr}\right\vert 
=  & - r_{n+1}(r_n\alpha +s_n \beta +t_n) + r_n(r_{n+1}\alpha +s_{n+1}\beta +t_{n+1}).
\cr}
$$
It follows that
$$
M(\uQ_n) \le h_{n+1}L(\uDelta_n) + h_nL(\uDelta_{n+1}) \le 2 h_{n+1} L_n =  2 h_n^{-v_n +v_n/w_n}.
$$
Bounding  from above the norm
$$
\| \uQ_n\| \le 2  \|\uDelta_n\| \|\uDelta_{n+1}\| \le  2 h_n h_{n+1} = 2 h_n^{1 +v_n/w_n}, 
$$
we find that
$$
M(\uQ_n) \le 2 (\| \uQ_n\|/2)^{-v_n(w_n-1)/(v_n  + w_n)}.
$$
For any $\epsilon >0$, we know that $w_n \ge w -\epsilon$, provided $n$ is large enough.
Selecting  an arbitrarily large index $n$ such that $v_n$ is arbitrarily close to the upper limit $v$, we 
obtain the lower bound $v'\ge v(w -1)/(v+w )$.

\medskip

The proof of the inequality  $v'\le (v-w+1)/w$ is quite similar, making now use of  the sequence  $(\uP_n)_{n\ge 1}$.
Define the non-zero integer triple $\uD_n$ by
$$
\uD_n = \uP_n\wedge \uP_{n+1} = \left( \left\vert \matrix{b_n & b_{n+1}\cr c_n & c_{n+1}\cr}\right\vert , 
- \left\vert \matrix{a_n & a_{n+1}\cr c_n   & c_{n+1}\cr}\right\vert , 
\left\vert \matrix{a_n & a_{n+1}\cr b_n & b_{n+1}\cr}\right\vert \right),
$$ 
so that $D_n$ is the line joining $P_n$ and $P_{n+1}$. Writing
$$
\displaylines{
 \left\vert \matrix{b_n & b_{n+1}\cr c_n & c_{n+1}\cr}\right\vert = 
 \left\vert \matrix{b_n -c_n\beta & b_{n+1}-c_{n+1}\beta \cr c_n & c_{n+1}\cr}\right\vert,
 \quad
 \left\vert \matrix{a_n & a_{n+1}\cr c_n   & c_{n+1}\cr}\right\vert
 = \left\vert \matrix{a_n -c_n \alpha & a_{n+1}-c_{n+1}\alpha \cr c_n   & c_{n+1}\cr}\right\vert , 
 \cr
 \left\vert \matrix{a_n & a_{n+1}\cr b_n & b_{n+1}\cr}\right\vert 
 =  \left\vert \matrix{a_n -c_n\alpha & a_{n+1}-c_{n+1}\alpha \cr b_n   & b_{n+1}\cr}\right\vert
+   \alpha \left\vert \matrix{c_n & c_{n+1}\cr b_n - c_n\beta    & b_{n+1}- c_{n+1}\beta \cr}\right\vert,
}
$$
and 
$$
 \left\vert \matrix{b_n & b_{n+1}\cr c_n & c_{n+1}\cr}\right\vert \alpha
 - \left\vert \matrix{a_n & a_{n+1}\cr c_n   & c_{n+1}\cr}\right\vert \beta 
 + \left\vert \matrix{a_n & a_{n+1}\cr b_n & b_{n+1}\cr}\right\vert 
 = \left\vert \matrix{c_n\alpha - a_n & c_{n+1}\alpha - a_{n+1}\cr c_n\beta - b_n & c_{n+1}\beta - b_{n+1}\cr}\right\vert,
 $$
 we obtain the upper bounds
 $$
 \displaylines{
 \| \uD_n\|     \le (1+ |\alpha|)\Big(  q_{n+1} M(\uP_n) +q_n M(\uP_{n+1})\Big) \le 2(1+|\alpha|) q_{n+1}^{1-w'_n},
 \cr
 L(\uD_n)   \le 2 M(\uP_n)M(\uP_{n+1}) \le 2 q_{n+1}^{-(v'_{n+1}+w'_n)},
 \cr}
 $$
 from which we deduce  the expected inequality 
 $$
 v\ge {v'+ w' \over 1 -w'} = v'w +w-1, 
 $$
taking into account  Jarn\'\i k's  relation $w'=(w-1)/w$.

\section
{3. The inverse problem.}

We have  to construct a point  $\Theta\in \bR^2$  for which the quadruple 
$\Omega(\Theta)$ takes a prescribed value $(v,v',w,w')$ as in the theorem. 
 We restrict in this part  to real numbers $w,w',v,v'$. The case of possibly infinite exponents
 is postponed to Section 7.  To this end,  we shall
 establish along  Sections 4--6 the following
 
 \proclaim
Proposition.
Let $w, \tau_0,\tau_1,\sigma$ be  positive real numbers satisfying the inequalities
$$
w \ge 2, \quad   \tau_1 \le 1 , \quad w \tau_0 \le \sigma \le \tau_0 +\tau_1.
$$
Then there exists $\Theta \in \bR^2$ such that
$$
\Omega(\Theta) = \Big(  {w-1+\tau_1\over \tau_0},{w-1\over \sigma},w,{w-1\over w}\Big).
$$

Let us show that the quadruples 
$$
(v,v', w, w') = \Big ({w-1+\tau_1\over \tau_0}, {w-1\over \sigma}, w,{w-1\over w}\Big)
$$ 
given  by the proposition are exactly those for which  the conditions 
$$
w\ge 2, \quad  w' = {w-1\over w} \and {v(w-1)\over v+w} \le v' \le {v-w+1\over w}
$$
of our theorem hold.

Let us fix $w\ge 2$. Observe first that, for given real numbers $\tau_0> 0$ and $0 < \tau_1 \le 1$,  the interval 
$$
w\tau_0 \le \sigma \le \tau_0 + \tau_1
$$
occurring in the proposition,  is non-empty
exactly when $(\tau_0,\tau_1)$ belongs to  the triangle ${\cal T}  \subset   \bR^2$  defined by the inequalities
$$
1 \ge \tau_1 \ge (w-1)\tau_0>0 . 
$$
 Fix now  $v \ge w(w-1)$. The necessity of this last assumption follows from Corollary 1.
 Then the intersection of $\cal T$ with the line of equation
 $$
 v\tau_0 = w-1 + \tau_1
 $$
in the plane $\bR^2$, is the segment whose extremities  are the points
$$
\Big({w-1\over v-w+1}, {(w-1)^2\over v-w+1}\Big)Ê\and \Big({w\over v},1\Big).
$$
The set of all admissible values $\sigma$, when the point $(\tau_0,\tau_1)$ ranges along this segment, coincides  with   the interval
$$
{w(w-1)\over v-w+1} \le \sigma \le {w\over v} +1 .
$$
  Therefore 
$v'=(w-1)/\sigma$ takes any assigned value in the interval
$$
{v(w-1)\over v+w} \le v' \le {v-w+1\over w}.
$$

In order to construct  a point $\Theta$ as in the proposition, it will be relevant
to assume  the  stronger conditions 
$$
w \ge 2, \quad  0 < \tau_0 < \tau_1 \le 1 , \quad w \tau_0 \le \sigma \le \tau_0 +\tau_1 ,\quad \sigma  < w-1 + \tau_0.
\leqno{(1)}
$$
   Notice that the assumptions of the proposition   imply  the slightly weaker inequalities 
$$
 0 < \tau_0 \le \tau_1 \le 1 \and \sigma \le w-1 + \tau_0.
 $$
  Hence   the additional  constraints in $(1)$ 
 exclude only the choices of parameters
$$
w =2, \quad \tau_0=\tau_1 ,\quad  \sigma = 2 \tau_0  \and   w=2, \quad \tau_1=1,\quad \sigma = 1 + \tau_0,
$$
which lead to  extremal quadruples of the form
$$
( v, {v-1\over 2},2,{1\over 2} ) \quad {\rm and } \quad ( v, {v\over v + 2},2,{1\over 2} )
$$
for some $v\ge 2$.  It turns out that Jarnik [\JarA--\JarC] has established    for any $v\ge 2$ 
the existence of points $\Theta$ for which
$$
(\om(\Theta)  , \om({}^t\Theta)) = (v,{v-1\over 2}) \and (\om(\Theta)  , \om({}^t\Theta)) = (v,{v\over v+ 2}).
$$
Then  we deduce from our refined transference inequalities that 
$$ \omc(\Theta)= 2 \and 
\omc({}^t\Theta) = 1/2.
$$
  We shall therefore assume  that (1) holds without any  loss of generality.

\section
{4. Constructing points and lines in the plane.}

We shall construct  in the next section a sequence  of points and a sequence of lines
which may be viewed as analogues of the sequences $(P_n)_{n\ge 1}$ and $(\Delta_n)_{n\ge 1}$  considered in  Section 2. 
To that aim, we establish  here  preliminary results.
 Lemma 1 provides us with
    families of rational points which  are close together and lie  on a given rational line. 
Next, we rephrase  dually our result  to obtain     families of close rational lines passing through  a given rational point.
As a main tool, we take again   standard arguments  arising from the theory of continued fractions.

Let us first introduce various notions of distances between  points and lines in the projective plane $\bP^2(\bR)$, and state 
some of their (easily proved) properties. It is convenient to view $\bR^2$ as a subset of $\bP^2(\bR)$
via the usual embedding $(x,y)\mapsto (x:y:1)$. With some abuse of notation, we shall  identify
a point in $\bR^2$ with its image in $\bP^2(\bR)$.  

For any pair of points $P$ and $P'$   
 in $\bP^2(\bR)$, with   homogeneous coordinates $\uP$ and $\uP'$,   denote by 
 $$
  d (P,P') =  {\| \uP \wedge \uP' \| \over  \| \uP\| \| \uP'\|}
   $$ 
   the so-called {\it projective distance} between $P$ and $P'$, which is obviously independent on the choice of $\uP$ and $\uP'$.
 The projective distance coincides inside the square $[-1/2,+1/2]^2$ with the distance associated to the norm of supremum. 
 In other words,  the formula
 $$
 d((x,y),(x',y')) = \max( | x-x'| , | y-y' | ) 
 $$
 holds whenever $\max(| x |,|y | ) \le 1/2$ and  $\max(| x' |,|y' | ) \le 1/2$. Moreover, for any $0\le R < 1$, the  projective ball 
 $$
 \Big\{ P\in \bP^2(\bR) \, ; \,\,  d(P,(0,0)) \le R\Big\},
 $$
 centered at the origin of $\bR^2$ with radius $R$,  is  equal to the square $[-R,+R]^2$. Note also that the triangle inequality
 $$
  d(P,P') - 2 d(P',P'') \le d(P,P'') \le d(P,P')+ 2 d(P',P'')
 $$
 holds  for any points $P,P',P''$ in $\bP^2(\bR)$ (see formula (5) of [\RoyC]).
   Now,  let  $\Delta$ be a  line in $\bP^2(\bR)$  with equation $rx+sy+tz=0$.  
 We set $\uDelta = (r,s,t)$
 and   define the (projective) distance $d(\Delta,\Delta')$
between two lines  $\Delta$ and $\Delta'$ by the formula
$$
d(\Delta , \Delta') = {\| \uDelta \wedge \uDelta'\|\over \| \uDelta\| \| \uDelta'\|}. 
$$
The distance $d(\Delta,\Delta')$  is again independent on the choice of the  triples $\uDelta$ and $\uDelta'$
respectively associated to $\Delta$ and $\Delta'$. 
Suppose that $\Delta$ intersects  $\Delta'$ inside the square  $[-1,+1]^2$. 
Then, denoting by   $ \langle \Delta , \Delta'\rangle$ the
acute angle  determined by the two lines in $\bR^2$, 
we have  the estimate
$$
{1\over 2}  \sin \langle \Delta , \Delta'\rangle \le  d(\Delta , \Delta')   \le  \,  2  \sin \langle \Delta , \Delta'\rangle.
$$
Finally,  we define the distance 
$d(P,\Delta)$,  between  a point $P$ with homogeneous coordinates $\uP= (x,y,z)$,    and  a line $\Delta$ with leading coefficients
$\uDelta =(r,s,t)$,    as  being   the quantity
$$
d(P,\Delta)  ={| r x  + s y +t z | \over \| \uP\|  \| \uDelta \|} .
$$
In the next sections, we shall make use  of the  formula
$$
 d(P,\Delta) = d(P,P')\, d(\Delta,\Delta'), \leqno{(2)}
$$
which is valid for any point  $P' $  of $\Delta$, distinct from $P$,
 and where $\Delta'$ stands for the line joining $P$ and $P'$.
This  equality, which follows from the formula for the  double 
vector product in $\bR^3$,   shows moreover that $d(P,\Delta)$  compares with  the
minimal projective distance between $P$ and the points of $\Delta$. 

  We call {\it normalized}  homogeneous coordinates of a rational point $P$ in $\bP^2(\bR)$, 
   any  triple $\uP =(a,b,c)$ of homogeneous coordinates of $P$, such that  $a,b,c$ are  coprime integers. 
 The triple $\uP$ is clearly defined up to a multiplicative factor $\pm 1$,  and we denote  by 
    $$
H(P)= \| \uP\| =\max(| a|, | b|, | c|)
$$ 
the usual height of the rational point $P$.  Note  that $H(P) =| c|$ when $P$ is located in the unit square $[-1,+1]^2$. 
Similarly,  we normalize the equation $rx+sy+tz =0$ of  any  rational projective line $\Delta$
by requiring that  $r,s,t$ are coprime integers, and we define its height $H(\Delta)$  as
being  the norm   $ \max(| r|, | s|, | t|)$. 
 
  \proclaim
 Lemma 1. Let $\Delta$ be a rational line in $\bP^2(\bR)$ with height $h$, 
 and let $P_0$
  be a rational point belonging to $\Delta$ with height $q_0$.  Let $\ell$ be a positive integer
 and let $q_1, \dots , q_\ell$ be a sequence of positive real numbers satisfying
 $$
 \quad q_{1} \ge 14 q_0 , \quad  q_0q_1 \ge  4 h \and q_{k+1}\ge 3 q_k  \quad (0\le k \le \ell -1).\leqno{(3)}
 $$
 There exist rational points $P_1, \dots , P_\ell $ located on $\Delta$,  such that   the estimates
 $$
{1\over 2} q_k \le H(P_k) \le 2  q_k   \quad ( 0 \le k \le \ell )
 $$
 and
 $$
  {1\over 32}   {h \over q_k q_{k+1}} \le d(P_{k}, P_{k'}) \le  16 {h\over q_k q_{k+1}} \quad (0 \le k < k' \le \ell )
 $$
 are verified. On the other hand, for any pair of distinct rational points $P$ and  $P'$ lying on $ \Delta$,  we have the lower bound
 $$
 d(P,P') \ge  {h\over H(P) H(P') }  .
 $$

 \proof
 Fix an equation  $rx+sy+tz=0$  of $\Delta$ whose coefficients $r,s,t$ are coprime integers, so that $h= \max(| r|, | s|, | t|)$.
 We  denote by $\Delta(\bZ)$ the additive group of integer triples $(a,b,c)$ for which   $ra + sb+ tc=0 $.
  Then  a rational   point $P$ lies on $\Delta$, if and only if  its normalized homogeneous coordinates $\uP$ belong  to $\Delta(\bZ)$.
  Thanks to Theorem  2  in [\BoVa], we may find  a basis $ \uA ,\uB$ of the $\bZ$-module  $\Delta(\bZ)$ such that 
$$
\|\uA\|\le \| \uB \|    \and    \|\uA\| \|\uB\| \le \sqrt{3}   h .
$$
An  integer triple $ m \uA + n\uB$  is primitive if and only if the coefficients $m$ and $n$ are relatively prime integers.
Note  that  $\uA\wedge \uB =\pm(r,s,t)$, so that $\| \uA\wedge\uB\| =h$.

 We first prove  {\it Liouville's inequality}, which is the last assertion of Lemma 1. 
  Suppose  that $P$ and $P'$ are distinct rational points located on  the  line $\Delta$. 
  Let $\uP$ and $\uP'$ be normalized homogeneous coordinates of $P$ and $P'$. 
  Then the vector  product $\uP\wedge \uP'$  is a non-zero integer multiple of $\uA\wedge \uB$.
  Therefore, we obtain the lower bound
  $$
  d(P,P')= 
  { \| \uP\wedge \uP' \| \over  \| \uP\| \|\uP'\| } \ge  {h \over H(P) H(P') } .
  $$
       
  Let $\uP_0$ be  normalized homogeneous coordinates of the point $P_0$. 
  Since $\uP_0$ belongs to $\Delta(\bZ)$, we may  write
$
\uP_0  = m \uA + n \uB
$
 for some  coprime integer coefficients $m$ and $ n$. Using Cramer's formula, we easily obtain the upper bounds 
$$
| m | \le  {2 q_0 \| \uB\| \over h} \le {2 \sqrt{3} q_0\over \| \uA \|} \and | n | \le {2 q_0 \| \uA\| \over h} \le {2 \sqrt{3} q_0\over \| \uB\|}.
$$
 Let    $e$ and $f$ be integers satisfying the equation
$
m f -  n e =  1, 
$
chosen so  that $f$ has minimal absolute value. Suppose first that $n$ is non-zero. Noting that $f$ is an  element of smallest absolute
value  in some coset modulo  $n$, we bound
$$
| f| \le {| n| \over 2} \le {\sqrt{3} q_0 \over \| \uB\|} \and
| e| \le { | m|| f | +1\over | n|} \le {| m_|\over 2} +1 \le {\sqrt{3} q_0\over \|\uA\|} +1.
$$
Thus 
$$
\| e \uA+ f \uB \|  \le 2\sqrt{3} q_0 +\| \uA \| \le 4\sqrt{3} q_0 \le q_1/2,
$$
since $ \| \uA \| \le \| \uB \| \le | n | \| \uB \| \le 2 \sqrt{3} q_0$. When $n=0$, we have  $\uA =\pm\uP_0$. Then $e=0, f=\pm 1$, and 
we   bound   again 
$$
\| e \uA+ f \uB \| =  \| \uB \| \le {\sqrt{3} h\over \| \uA\|}\ = {\sqrt{3} h\over q_0} \le  q_1/2.
$$

  We are now able to  construct  the sequence of points $P_1,\dots ,P_\ell$.  Define 
    $$
 g_1=  \left\lceil {q_1  \over q_0 } \right\rceil  \and  \uP_1 =  g_1 \uP_0 +e \uA + f \uB .
    $$
 The integer triple $\uP_1$ is primitive. Let $P_1$ be  the rational point in $\bP^2(\bR)$ with  homogeneous coordinates $\uP_1$.
 Its height $H(P_1)$ is therefore  equal to $ \| \uP_1\|$,  and satisfies the required estimate  
  $$
   q_1/2 \le q_1 - \| e\uA + f\uB \| \le  \| \uP_1\| =H(P_1)  \le q_1+ q_0 + \| e\uA + f\uB \| \le 2 q_1.
   $$  
 Next, when $\ell \ge 2$, we define recursively a sequence of primitive 
  integer triples $\uP_2, \dots , \uP_\ell$ by the relations 
    $$
\uP_{k} = g_k \uP_{k-1} + \uP_{k-2}, \quad (2\le k \le \ell) ,
$$
  where we have set  
  $$
  g_k = \left\lceil {q_k  \over \| \uP_{k-1}\| } \right\rceil \ge 1 . 
  $$
Let  $P_k$ be the rational point with homogeneous coordinates $\uP_k$. 
  Arguing by induction on $k$, we obtain similarly the  estimate of height
$$
{1\over 2} q_k  \le q_k - \| P_{k-2}\| \le H(P_k) = \| \uP_k\| \le q_k +\| P_{k-1}\|  + \| P_{k-2}\|  \le 2 q_k.\leqno{(4)}
$$

It remains to estimate the distances between the points $P_k$. Let us write
$$
\uP_k= u_k \uP_0 + v_k \uP_1 
$$
for  integer coefficients $u_k, v_k$ satisfying the usual recurrence relations	
$$
u_k = g_k u_{k-1} + u_{k-2}, \quad (u_0 = 1, u_1 = 0) \and
v_k = g_k v_{k-1} + v_{k-2}, \quad (v_0 = 0, v_1 = 1),
$$
occurring in  the theory of continued fractions.  We  therefore have the formula 
$$
{u_k\over  v_k} = [0 ; g_2,\dots , g_k ], \quad (1 \le k\le \ell).
$$
Observe  next that the estimates of norms
$$
\eqalign{
{1\over 2} v_k \| \uP_1\|  \le  v_k   & \left( \| \uP_1\| - [0;g_2, \dots ,g_k] \|\uP_0\| \right)=   v_k \| \uP_1\| - u_k \| \uP_0\|  
 \le \| \uP_k\|  \le 
 \cr
 &     u_k \| \uP_0\|  + v_k \| \uP_1\|  =  v_k\left( [0;g_2, \dots ,g_k] \|\uP_0\| +  \| \uP_1\|\right)
\le 2 v_k \|\uP_1\| 
\cr}
 \leqno {(5)}
$$
hold as well for $1\le k \le \ell$. Now for any $0 \le k < k' \le \ell$, we have 
$$
d(P_k,P_{k'}) = { \| \uP_k\wedge \uP_{k'}\| \over \| \uP_k\| \| \uP_{k'}\| }
 ={ | u_kv_{k'} - u_{k'}v_k| h \over \| \uP_k\| \| \uP_{k'}\| },
$$
since 
$$
\displaylines{
\qquad \uP_k\wedge \uP_{k'} = (u_kv_{k'} - u_{k'}v_k) \uP_0\wedge \uP_1 \hfill
\cr
\hfill 
= (u_kv_{k'} - u_{k'}v_k)(mf-ne)\uA\wedge \uB = \pm(u_kv_{k'} - u_{k'}v_k) (r,s,l).\qquad
\cr}
$$
By a standard result on continued fractions, we know that
$$
{1\over 2 v_{k+1}} \le \Big | v_k{u_{k'}\over v_{k'}} -u_k \Big | \le {1\over  v_{k+1}} .
$$
 It follows that 
$$
{1\over 2}  {v_{k'} h \over v_{k+1}  \|\uP_k\| \| \uP_{k'}\| }  \le d(P_k,P_{k'}) \le   {v_{k'} h \over v_{k+1}  \|\uP_k\| \|\uP_{k'}\| } . 
$$
The required estimate for $d(P_k,P_{k'})$  follows  from  (4) and (5). 
\cqfd

\bigskip 

We state now a  dual version  of Lemma 1, in which the roles of lines and points are exchanged.

\proclaim
Lemma 2. 
Let  $\Delta_0 $ be a rational line with    height $h_0$, and let  $P$ be a rational point lying on $\Delta_0$,  with height $q$. 
 Let $\ell$ be a positive integer
 and let $h_1, \dots , h_\ell$ be a sequence of positive real numbers satisfying
 $$
 \quad h_{1} \ge 14 h_0 , \quad  h_0h_1 \ge  4 q \and h_{k+1}\ge 3 h_k  \quad (0\le k \le \ell -1).\leqno{(6)}
 $$
 There exist rational lines $\Delta_1, \dots , \Delta_\ell $ passing through the point $P$,  such that  the estimates of height
 $$
 {1\over 2} h_k \le H(\Delta_k) \le 2  h_k   \quad (0 \le k \le \ell )
 $$
 and of distance
 $$
  {1\over 32}   {q\over h_k h_{k+1}} \le d(\Delta_{k}, \Delta_{k'}) \le 16 {q\over h_k h_{k+1}} \quad (0 \le k < k' \le \ell )
 $$
 are verified.  On the other hand, for any pair of distinct rational lines $\Delta$ and  $\Delta'$ containing $ P$,  we have the lower bound
 $$
 d(\Delta,\Delta') \ge {q\over H(\Delta) H(\Delta') }  .
 $$

\proof
It is completely parallel to the proof of Lemma 1. The formalism of the proof remains exactly the same. We omit the details. \cqfd

\section
{5. The basic  construction.}
Remind the stronger assumptions (1) relating  the data $w,\tau_0,\tau_1,\sigma$. 
Observing that $0 < \tau_0 < \tau_1 \le 1$, we put $\ell =1$ if $\tau_1=1$,  and  pick otherwise an integer $\ell \ge 2$ and  
 an increasing sequence of  real numbers $\tau_2, \dots , \tau_\ell$,  such that
$$
0< \tau_0 < \tau_1 < \dots < \tau_\ell =1  \and
 {w-1+\tau_{k+1} \over \tau_k}  \le {w-1+\tau_1\over \tau_0} , \quad (0 \le k \le \ell -1). \leqno{(7)}
$$
As an example, we may choose $\tau_k = \min  ( 1, \tau_0 + k(\tau_1 -\tau_0))$ for $k= 1, \dots, \ell$, where $\ell$
is the smallest  integer such that $ \tau_0 + \ell( \tau_1 -\tau_0) \ge 1$. 
Note that, in any case, the sequence $(\tau_k)_{0\le k\le \ell}$ increases and ends at  $\tau_\ell= 1$. 
Set now
$$
\sigma_0 = \sigma \and \sigma_1 = w.
$$
It follows from $(1)$ that $0< \sigma_0 < \sigma_1 \le \sigma/\tau_0$. We extend similarly this second sequence  into an (eventually longer)
increasing sequence 
$$
0 < \sigma_0   <    \dots < \sigma_{\ell'} = {\sigma / \tau_0} , \leqno {(8)}
$$
for some integer $\ell' \ge 1$  selected  so  that the growth conditions
$$
 {\sigma_{k+1}-1\over \sigma_k} \le {w-1\over \sigma}, \quad (0 \le k \le \ell' -1) \leqno {(9)}
$$
hold. The constraints (8) and (9)  can be simultaneously fulfilled by taking bounded ratios
$$
1 < {\sigma_{k+1}\over \sigma_k} \le {w-1+\tau_0\over \sigma} \quad {\rm for}\quad 1 \le k\le \ell' -1.
$$

Next,   we introduce two  sequences of positive real numbers  $(h_{n,k})_{n\ge 1, 0\le k \le \ell }$ and $(q_{n,k})_{n\ge 1, 0\le k \le \ell'}$
in the following way.
For simplicity, set $h_n= h_{n,0}$. We start with any   large 
initial value $h_1$,  and  define  inductively $h_{n,k}$ and $q_{n,k}$ thanks to the recurrence relations
$$
 h_{n+1} =   h_n^{1/\tau_0}, \quad 
h_{n,k} =  h_{n+1}^{\tau_k},  (0\le k\le \ell) 
\and 
  q_{n,k} =  h_{n+1}^{\sigma_k}/16,  \,\,(0\le k \le \ell').
\leqno{(10)}
$$
Taking into account the equalities $ \tau_\ell=1$ and $\sigma_{\ell'}= \sigma/\tau_0$,
 observe  that the branching equations
 $$
h_{n,\ell}= h_{n+1,0}  \and  q_{n,\ell'}= q_{n+1,0} 
$$
 hold for any $n\ge 1$. With the exception of  these  equalities, the sequences  $(h_{n,k})$ and $(q_{n,k})$,
where the indices $(n,k)$ have been lexicographically ordered, are strictly increasing,  since so are the sequences of exponents
$(\tau_k)$ and $(\sigma_k)$. Notice as well that 
both sequences $(h_{n,k})$ and $(q_{n,k})$ increase at least as  a double exponential.

For further use, let us quote the estimates 
 $$
 {h_{n+1}\over h_n} < {q_{n,1}\over q_{n,0} } 
\and  
{ h_{n,k} h_{n+2}\over q_{n+1,0}q_{n+1,1}}  = o\Big({h_{n+1} \over h_{n,k+1} q_{n,1}}\Big) , \quad (0\le k\le \ell -1), 
\leqno{(11)}
$$
which  follow from (1) and $(10)$.  In order to check the second part of $(11)$, write both ratios in terms of $h_{n+1}$, and observe that
$$
\tau_k +{1-\sigma - w\over \tau_0} < 1 +{ 1-w\tau_0 -w\over \tau_0} = -(w-1)(1+{1\over \tau_0}) < -2(w-1) \le 1-\tau_{k+1}-w.
$$

\proclaim
Lemma 3.
There exists   a sequence of rational lines $(\Delta_{n,k})_{n\ge 1, 0\le k \le \ell}$
and  a sequence of rational points  $(P_{n,k})_{n\ge 1, 0\le k\le \ell'}$  satisfying
for any $n\ge 1$  the compatibility relations
$$
\Delta_{n,\ell}= \Delta_{n+1,0} , \quad  P_{n,\ell'}= P_{n+1,0},
$$
and  the following properties. 
The points $P_{n,0},\dots , P_{n,\ell'}$ are pairwise distinct
and lie on the line $\Delta_{n+1,0}$. The lines $\Delta_{n,0},\dots , \Delta_{n,\ell}$ are pairwise distinct and pass  through
the point $P_{n,0}$. Moreover
the estimates  \footnote{{\rm(*)}}{{\rm 
The implicit constants involved in the  forthcoming symbols $\gg$ and $\ll$ are absolute. 
Their computation is however useless for our purpose.
We frequently write $F(n)\gg\ll G(n)$ to 
signify that $F(n)\gg G(n)$ and $F(n)\ll G(n)$ for all sufficently large $n$.} 
}
of distances
$$
\eqalign{
 d( (0,0), P_{1,0} )  & \gg\ll {h_1\over q_{1,0}},
\cr
 d(P_{n,k},P_{n,k'})  & \gg\ll {h_{n+1}\over q_{n,k}q_{n,k+1}} , \quad (0\le k <k' \le \ell' ),
\cr
  d(\Delta_{n,k}, \Delta_{n,k'})  &\gg\ll { q_{n,0}\over h_{n,k}h_{n,k+1}}, \quad (0 \le k <k'\le \ell ),
\cr}\leqno{(12)}
$$
are satisfied,  as well as the estimates  of heights
$$
{1\over 2}q_{n,k} \le H(P_{n,k})  \le 2 q_{n,k}  \and {1\over 2}h_{n,k} \le  H(\Delta_{n,k}) \le 2 h_{n,k} . \leqno{(13)}
$$

\proof 
We carry out simultaneously the construction of  both sequences $\Delta_{n,k}$ and $P_{n,k}$ by successive steps.

 Start (for instance) by defining 
 $\Delta_{1,0}$   as  the line with  equation $ \lceil h_1 \rceil x - y =0$, and by choosing  the point 
 $P_{1,0}= (1/\lceil q_{1,0}\rceil , \lceil h_1\rceil /\lceil q_{1,0}\rceil )$  on $\Delta_{1,0}\cap \bR^2$.
 Observe that 
  $$
  q_{1,0} = h_2^{\sigma}/16 \ge   h_1^w /16
  $$  
  is much bigger than $h_1$,  when $h_1$ is large. 
 Then, the required estimates 
 $$
 \displaylines{
 d( (0,0), P_{1,0} )  \gg\ll { h_1 \over  q_{1,0}} \, , 
  \cr
 {1\over 2} h_1 \le H(\Delta_{1,0})= \lceil h_1\rceil  \le 2 h_1  \and   
 {1\over 2} q_{1,0}  \le H(P_{1,0})= \lceil q_{1,0}\rceil  \le 2 q_{1,0} \, ,
\cr}
 $$
 clearly hold for sufficiently large initial values $h_1$. Note  also that $P_{1,0}$ may be taken arbitrarily close to the origin $(0,0)$
 of $\bR^2$, provided $h_1$ is large enough. 
 
  Suppose now that $P_{1,0},P_{1,1}, \dots, P_{n,0}$ and $\Delta_{1, 0},\Delta_{1,1},\dots , \Delta_{n,0}$  have already been selected
 for some $n\ge 1$.
 We first use  Lemma 2, applied  to the point $P_{n,0}$ lying on the line $\Delta_{n,0}$,  and to the sequence $h_{n,1}, \dots , h_{n,\ell}$.
 The main assumption $H(\Delta_{n,0})h_{n,1}\ge 4H(P_{n,0})$ occurring in (6),  follows from (10, 13) and from the inequality
 $\sigma \le \tau_0 +\tau_1$ in (1). Therefore we may find rational lines $\Delta_{n,1}, \dots , \Delta_{n,\ell}=\Delta_{n+1,0}$ 
 passing through $P_{n,0}$, for which the third estimate of 
  (12) and the second one of (13) are verified.  Next, starting with   the point $P_{n,0}\in \Delta_{n+1,0}$,
 we apply Lemma 1 to the sequence $q_{n,1},\dots ,q_{n,\ell'}$.    We  obtain rational points $P_{n,1}, \dots , P_{n,\ell'}= P_{n+1,0}$
 belonging to $\Delta_{n+1,0}$ and  satisfying $(12)$ and $(13)$. Notice that the condition
 $H(P_{n,0})q_{n,1} \ge 4 H(\Delta_{n+1,0})$ occurring in the assumptions
 (3) of Lemma 1 is easily  fulfilled, since 
 $$
 q_{n,1} =  h_{n+1}^{\sigma_1}/16 \ge  h_{n+1}^2/16.
 $$
 The two sequences have thus been extended up to the rank $(n+1,0)$.
 \cqfd

   Let us show that the sequence of points $(P_{n,k})_{n\ge 1,0\le k\le \ell'-1}$ 
   furnished by Lemma 3, is a Cauchy sequence in $\bP^2(\bR)$.
 Observe  first that  the sequence $(h_{n+1}/q_{n,k}q_{n,k+1})_{n\ge 1,0\le k\le \ell'-1}$, which occurs  in (12), 
 is  decreasing when  the indices  $(n,k)$ are lexicographically ordered. The only non-obvious inequality
 $$
 {h_n\over q_{n-1,\ell' -1}q_{n-1,\ell'} } > {h_{n+1}\over q_{n,0}q_{n,1}}
 $$
 follows from the first part  of (11). Moreover this sequence tends  to $0$, 
  much more quickly than any geometrical sequence  with ratio $<1$.  
  Take any  index $(n,k)$  smaller  than  $(n',k')$ for that lexicographic order. 
 Combining the triangle inequality 
 with the upper bounds (12), we find 
 $$
 \eqalign{
 d(P_{n,k},P_{n',k'}) \le  & \sum_{ (n,k) \le (\nu ,\kappa)< (n',k')}2^{{\rm rk}(\nu ,\kappa)} d( P_{\nu ,\kappa},P_{\nu ,\kappa +1})
 \cr 
  \ll & \sum_{ (n,k) \le (\nu ,\kappa)< (n',k')}2^{{\rm rk}(\nu ,\kappa)} { h_{\nu  +1} \over q_{\nu , \kappa }q_{\nu ,\kappa +1}}
  \ll 
 { h_{n+1}\over q_{n,k}q_{n,k+1}},
 \cr }
 $$
where ${\rm rk}(\nu, \kappa)$ denotes the rank of   $(\nu ,\kappa)$ in the  ordered sequence $(n,k)< \dots <(n',k')$
of  all indices between $(n,k)$ and $(n',k')$,
 starting with the initial value ${\rm rk}(n,k)=0$. 
 
  Let $\Theta$ be the limit of the  Cauchy sequence $(P_{n,k})$. 
  The same argument  as above  yields the  estimates 
$$
d(P_{n,k}, \Theta)  \gg\ll {h_{n+1}\over q_{n,k}q_{n, k+1}}, \quad (0 \le k \le \ell' -1)  ,\leqno{(14)}
$$
and
$$
d(P_{n,0},P_{n+1,0}) \gg\ll {h_{n+1}\over q_{n,0} q_{n,1}}.
\leqno{(15)}
$$
Taking moreover $h_1$ large enough, we may  assume that $d((0,0), P_{n,k}) \le 1/4$ for any index $n\ge 1, 0\le k \le \ell' -1$, 
so that all points $P_{n,k}$ lie in the square $[-1/4,+1/4]^2$. Then $\Theta$ obviously belongs to $[-1/4,+1/4]^2$. 

Put now $\Theta =(\alpha,\beta)$ and recall  the notations
 $$
L(\uX) = | x\alpha + y\beta +z|, \quad  M(\uX) = \max (|z\alpha -x|,|z\beta -y|) , \quad \uX =(x,y,z),
$$
introduced in Section 2.   Let
$$
\eqalign{
\uP_{n,k} = & (a_{n,k}, b_{n,k}, c_{n,k}) 
\quad {\rm with} \quad q_{n,k}/2 \le | c_{n,k} | \le 2 q_{n,k}
\and \gcd (a_{n,k}, b_{n,k},c_{n,k})=1,
\cr
\uDelta_{n,k} = & (r_{n,k}, s_{n,k}, t_{n,k}) 
\quad {\rm with} \quad h_{n,k}/2  \le \|  \uDelta_{n,k} \| \le 2 h_{n,k}
\and \gcd (r_{n,k}, s_{n,k},t_{n,k})=1,
\cr}
$$
be  normalized integer triples respectively  associated  to   the rational point $P_{n,k}$ and to the rational  line $\Delta_{n,k}$.
Recall also that the projective distance $d$ coincides inside   the square $[-1/2,+1/2]^2$ with the distance associated to the norm of supremum.  
The estimate  of distance (14) is  therefore equivalent to 
$$
M(\uP_{n,k}) = | c_{n,k} | d(\Theta, P_{n,k}) \gg\ll  {h_{n+1}\over q_{n,k+1}} \gg\ll h_{n+1}^{1-\sigma_{k+1}}, \quad (0\le k \le \ell' -1).
\leqno{(16)}
$$
Observe next that  the point 
$
P_{n+1,0}=(a_{n+1,0}/c_{n+1,0}, b_{n+1,0}/c_{n+1,0})
$
 belongs to the line $\Delta_{n+1,0}$ which intersects  $\Delta_{n,k}$
at the point $P_{n,0}$. Employing now the formula  (2) to estimate the distance $d(P_{n+1,0},\Delta_{n,k})$, we 
deduce from (12, 13, 15)  that 
$$
\displaylines{ \quad
 |r_{n,k}{a_{n+1,0}\over c_{n+1,0} }  +s_{n,k} {b_{n+1,0}\over c_{n+1,0} }  + t_{n,k}|  \, \gg\ll \, 
  \| \uDelta_{n,k} \|  \,  d(\Delta_{n,k}, \Delta_{n+1,0}) \, d(P_{n,0}, P_{n+1,0})
\hfill\cr\hfill
\gg\ll {h_{n+1}\over h_{n,k+1}q_{n,1}} \gg\ll
h_{n+1}^{1-w-\tau_{k+1}}.
\quad\cr}
$$
Using  (14) and the second part of $(11)$, we may replace in the above inequalities $a_{n+1,0}/c_{n+1,0}$ 
and $ b_{n+1,0}/c_{n+1,0}$ by their
 limits  $\alpha$ and $\beta$. We  therefore obtain  the estimate
$$
L(\uDelta_{n,k}) \gg\ll  {h_{n+1}\over h_{n,k+1}q_{n,1}} \gg\ll
h_{n+1}^{1-w-\tau_{k+1}}, \quad (0\le k \le \ell -1).
\leqno{(17)}
$$

At the present stage, we have constructed two sequences of integer triples $\uP_{n,k}$ and $\uDelta_{n,k}$,
  which provide
good approximations to $\Theta$ with respect to the functions $M$ and $L$. 
Since the norm of $\uP_{n,k}$ (resp. $\uDelta_{n,k}$) compares  to $h_{n+1}^{\sigma_k}$ (resp.  $h_{n+1}^{\tau_k}$), 
 the upper bounds given  by  (16--17) yield the lower bounds
$$
\eqalign{
\om(\Theta) \ge & \max_{0\le k \le \ell -1}\left({ w-1 + \tau_{k+1}\over \tau_k}\right) = {w-1+ \tau_1\over \tau_0},
\cr
\om({}^t \Theta) \ge &  \max_{0\le k\le \ell' -1}\left({ \sigma_{k+1} -1\over \sigma_{k}}\right) = {w-1\over \sigma}.
\cr
\omc(\Theta) \ge & \min_{0\le k\le \ell -1}\left( {w-1 +\tau_{k+1}\over \tau_{k+1}} \right)= w,
\cr
\omc({}^t \Theta) \ge &  \min_{0\le k\le \ell' -1}\left({ \sigma_{k+1} -1\over \sigma_{k+1}}\right) = {w-1\over w}.
\cr}
\leqno {(18)}
$$
Notice that the two first  equalities on the right-hand side of $(18)$ follow  
 from (7--9). 

 It turns out that the  lower bounds (18) are actually equalities,   as  we shall prove  in the  next  section.

\bigskip
\section
{6. Upper bounds.}

In order to bound from above the exponents $\om({}^t \Theta)$ and  $\omc({}^t \Theta)$,  (resp. $\om(\Theta)$ and  $\omc(\Theta)$), 
we establish that the rational points,  (resp. the rational  lines),  which well approximate the point $\Theta$ belong necessarily
to the set of points $P_{n,k}$, (resp. the set of lines $\Delta_{n,k}$),  previously considered. 
That is the underlying principle for the proof of the next two lemmas. 

\proclaim
Lemma 4. For any non-zero integer triple $\uP$  which is not proportional to some  triple $\uP_{n,k}$,
and   having  sufficiently large norm $\| \uP\|$,    we have the lower bound 
\footnote{{\rm (*)}}{{\rm 
 In this section, the  constants involved in the symbols $\gg$ and $\ll$ may  possibly depend
 upon the data $w,\tau_0,\tau_1, \sigma$.}}
$$
M(\uP) \gg \| \uP \|^{-\lambda}\quad {\rm with } \quad \lambda = 
\max\Big( {1\over w-1 +\tau_0}, {\sigma -\tau_0\over \sigma }\Big) .
$$
 There exists a positive real number $\epsilon$ such that for  any sufficiently large integer $n$
 and for any non-zero integer triple $\uP$ with norm $\|\uP\| \le \epsilon q_{n,1}$, we have the uniform lower bound
$$
M(\uP)\gg \epsilon q_{n,1}^{-(w-1)/w}.
$$

\proof
Let us first observe that the sequence
$$
\gamma_{(n,k)} =  q_{n,1}{  h_{n,k}\over h_{n+1}}, \quad ( n\ge 1 , 0 \le k \le \ell),
$$
indexed by the couples of integers $(n,k)$ lexicographically ordered,   increases
(the inequality $\gamma_{(n+1,0)} > \gamma_{(n,\ell)}$  follows  from the first part of (11)), and tends to infinity.
We denote by $(n,k)+1$ the successor of $(n,k)$ for the lexicographic order.  Then the estimate  (17) may be written as
$$
 L(\uDelta_{n,k})\gg\ll  {1\over \gamma_{(n,k+1)}} = {1\over \gamma_{(n,k)+1}} , \quad (0\le k \le \ell -1).
$$

Write $\uP =(a,b,c)$ and put $q=\|\uP\|$. Suppose first that the point $P$ with homogeneous coordinates $\uP$ 
lies ouside the square $[-1/2,+1/2]^2$. 
Since $\Theta$ belongs to $[-1/4,+1/4]^2$, we bound from below
$$
M(\uP) = | c | \max\left( | \alpha -{a\over c}|, | \beta - {b\over c}|\right) \ge {1\over 4} ,
$$
when  $c$ is non-zero,  and $M(P) \ge 1$ if  $c=0$. 
We shall  therefore assume  that $P$ belongs to the square $[-1/2,+1/2]^2$, so that
$q= | c|$. The identity
$$
r_{n,k}a +s_{n,k}b + t_{n,k}c  = c(r_{n,k}\alpha +s_{n,k}\beta + t_{n,k}) -r_{n,k}(c\alpha -a) -s_{n,k}(c \beta -b)
$$
yields for any index $(n,k)$ the upper bound
$$
 |r_{n,k}a +s_{n,k}b + t_{n,k}c| \le   q L(\uDelta_{n,k})  +2 h_{n,k} M(\uP)
 \ll   {q \over \gamma_{(n,k)+1}} + h_{n,k} M(\uP).
\leqno {(19)}
$$

Let $\epsilon$ be a positive real number.  Assuming $q$ large enough,  we 
define $(n,k)$ as the unique index for which
$$
\epsilon \gamma_{(n,k)} < q \le \epsilon \gamma_{(n,k)+1}.
$$
We use a different argumentation whether $0\le k \le \ell -1$ or $k=\ell$.

Suppose first that $k\le \ell -1$. Replacing in (19) the index $k$ by $k+1 $, we obtain the two upper bounds
$$
\eqalign{
& |r_{n,k}a +s_{n,k}b + t_{n,k}c| \ll \epsilon + h_{n,k} M(\uP)
\cr 
&  |r_{n,k+1}a +s_{n,k+1}b + t_{n,k+1}c| \ll \epsilon + h_{n,k+1} M(\uP) .
\cr}
 \leqno {(20)}
 $$
If we suppose that $M(\uP)\le \epsilon h_{n,k+1}^{-1}$ and $\epsilon$ small enough, 
 the left-hand sides of both inequalities  (20) must vanish, since these are integers. It follows that 
$$
P = \Delta_{n,k}\cap \Delta_{n,k+1} = P_{n,0},
$$
which is impossible since we have assumed that $P$ differs from all  points $P_{n,k}$. 
Therefore 
$$
M(\uP) > \epsilon  h_{n,k+1}^{-1} = \epsilon h_{n+1}^{-\tau_{k+1}}
\and q  > \epsilon \gamma_{(n,k)} =  {\epsilon\over 16} h_{n+1}^{w-1 +\tau_k},
$$
so that 
$$
M(\uP) \gg \epsilon (q/\epsilon)^{-\tau_{k+1}/(w-1+\tau_k)} \gg \epsilon (q/\epsilon)^{-\lambda}.
$$

We consider now the case $k=\ell$. Then $q$ belongs  to the interval 
$$
\epsilon \gamma_{(n,\ell)} = \epsilon q_{n,1} < q \le \epsilon q_{n+1,1}{h_{n+1}\over h_{n+2}}
=\epsilon \gamma_{(n+1,0)}.
$$
In this situation, (19)  yields the upper bound
$$
 |r_{n+1,0}a +s_{n+1,0}b + t_{n+1,0} c| \ll  {q\over \gamma_{(n+1,1)}} + h_{n+1} M(\uP) \ll \epsilon +h_{n+1}M(\uP).
 $$
 When the point $P$ does not lie on the line $\Delta_{n+1,0}$, the  sum  $r_{n+1,0}a +s_{n+1,0}b + t_{n+1,0} c$
  is a non-zero integer.
 Thus, if $\epsilon$ is small enough,  we obtain in this case the stronger lower bound
 $$
 M(\uP) \gg h_{n+1}^{-1} = (16 q_{n,1})^{-1/w} \gg (q/\epsilon)^{-\lambda}.
 $$
 It remains to deal with points $P$ located on $\Delta_{n+1,0}$. Since $\epsilon q_{n,1} < q$,  define  $k'$
 as  the largest positive integer $k\le \ell'$ such that
 $
 \epsilon q_{n,k} <q.
 $
 Therefore  
 $
 \epsilon q_{n,k'}< q\le \epsilon q_{n,k'+1}
 $
 when $1\le k' \le  \ell' -1$,  and 
 $
 \epsilon q_{n+1,0} < q \le \epsilon q_{n+1,1}h_{n+1}/ h_{n+2}
 $
 when $k'=\ell'$. Now, Liouville inequality (on the line $\Delta_{n+1,0}$) provides us with the lower bound
 $$
 d(P,P_{n,k'}) \gg   {h_{n+1}\over q \, q_{n,k'}} \gg \epsilon^{-1}{h_{n+1}\over q_{n,k'}q_{n,k'+1}} 
  $$
 when $1 \le k' \le \ell' -1$, or
 $$
 d(P,P_{n,\ell'}) \gg  {h_{n+1}\over q \, q_{n+1,0}} \gg \epsilon^{-1} {h_{n+2}\over  q_{n+1,0} \, q_{n+1,1}}  
 $$
 when  $k'= \ell'$.  On the other hand, (14) gives the upper bounds
 $$
 d(P_{n,k'} , \Theta) \ll {h_{n+1}\over q_{n,k'}q_{n,k'+1}}, \quad  (1 \le k' \le \ell' -1)
 \and
 d(P_{n,\ell'}, \Theta)  \ll {h_{n+2}\over  q_{n+1,0} \, q_{n+1,1}}   .
 $$
 In both cases, the triangle inequality  shows that
 $$
 d(P,\Theta) \gg {h_{n +1}\over q \, q_{n,k'} } ,
 $$
 provided $\epsilon$  is small enough. It follows  that 
 $$
 M(\uP) = q \, d(P,\Theta) \gg {h_{n+1} \over q_{n,k'}} = 16 h_{n+1}^{1- \sigma_{k'}}\gg (q/\epsilon)^{-(\sigma_{k'} -1)/\sigma_{k'}}
 \gg (q/\epsilon)^{-\lambda},
 $$
 noting that the exponent $1-\sigma_{k'}$ is negative, since $k'\ge 1$ and $\sigma_{k'} \ge \sigma_1 =w\ge 2$. 
 Fixing finally $\epsilon$ small enough  so that the previous estimates are valid, we have thus proved the first assertion of Lemma 4.  
 
 As for the second part of Lemma 4, we take again the same argumentation
in a simpler way.  Observe that 
$$
 \epsilon q_{n,1} = \epsilon \gamma_{(n, \ell)} = {\epsilon\over 16} h_{n+1}^{w}.
$$
 Let $\uP=(a,b,c)$ be any  non-zero integer triple with norm $\|Ê\uP\| \le \epsilon q_{n,1}$.
Now  (19)  gives 
$$
\max_{k\in \{\ell -1, \ell\} }  |r_{n,k}a +s_{n,k}b + t_{n,k}c|
\ll \epsilon +h_{n+1} M(\uP).
$$
If $M(\uP) \le \epsilon h_{n+1}^{-1}$  and   $\epsilon$ small enough,  the left-hand side of the above inequality vanishes, 
and we find that
$$
P = \Delta_{n,\ell -1}\cap \Delta_{n, \ell} = P_{n,0}.
$$
  Then by (16)
$$
M(\uP) \ge M(\uP_{n,0}) \gg h_{n+1}^{1-w} \gg   q_{n,1}^{-(w-1)/w}.
$$
  Otherwise 
$$
M(\uP) > \epsilon h_{n+1}^{-1} \gg \epsilon q_{n,1}^{-1/w} .
$$
Therefore  the lower bound $M(\uP) \gg \epsilon q_{n,1}^{-(w-1)/w}$ holds for any non-zero integer triple
$\uP$ with norm  $\| \uP\|  \le \epsilon q_{n,1}$. 
\cqfd

\bigskip

The next result may be viewed as a dual version of Lemma 4.

\proclaim
Lemma 5. For any non-zero integer triple   $\uDelta$  whose norm $ \| \uDelta \|$ is large enough,  and
which is  not  proportional to some triple $\uDelta_{n,k}$, we have the lower bound
 $$
L(\uDelta)  \gg \| \uDelta \|^{-\mu} \quad {\rm with} \quad 
\mu =\max\Big( {\sigma\over (w-1) \tau_0}, {w-1+\tau_0\over \tau_0}\Big).
$$
There exists a positive real number $\epsilon$ such that for  any sufficiently large integer $n$ and for any non-zero integer triple  $\uDelta$
with norm $\le \epsilon h_n $, we have the uniform lower bound
$$
L(\uDelta) \gg  h_n^{-w}.
$$

\proof
We take again  the same arguments as in Lemma 4, exchanging the roles of lines and points. 
Set now (note that  $k\ge 1$ here)
$$
\delta_{(n,k)} =  {q_{n-1,k} \over  h_{n} }, \quad (n\ge 2, 1 \le k \le \ell').
$$
The sequence $\delta_{(n,k)}$,
indexed by the couples of integers $(n,k)$ with $1\le k \le \ell'$ in lexicographical order,   increases   and tends to infinity.
We denote again by $(n,k)+1$ the successor of $(n,k)$ relatively to the lexicographic  order.  Notice that  (16) may  actually be written 
in the form
$$
M(\uP_{n-1,k} ) \gg\ll  {1\over \delta_{(n,k)+1}} , \quad (1\le k \le \ell').
$$

Let $\epsilon$ be a positive real number which will be selected later sufficiently small.
Write $\uDelta=(r,s,t)$ and put $h=\| \uDelta \|$. Assuming $h$ large enough,
there exists a unique  index  $(n,k)$ such that
$$
\epsilon \delta_{(n,k)} < h \le \epsilon \delta_{(n,k)+1}.
$$
A similar  splitting of cases occurs as in Lemma 4. 

Suppose first  that $1 \le k\le \ell' -1$. Then we bound from above
$$
\eqalign
{
&  | r a_{n-1,k} + s b_{n-1,k} +   t c_{n-1,k} | \le   |  c_{n-1,k} |  L(\uDelta)  +2 h M(\uP_{n-1,k})
 \ll   q_{n-1,k} L(\uDelta)  + \epsilon ,
 \cr
&  | r a_{n-1,k+1} + s b_{n-1,k+1} +  t c_{n-1,k+1} | 
 \ll     q_{n-1,k+1}  L(\uDelta)  + \epsilon .
\cr
}
\leqno {(21)}
 $$
If we suppose that $L(\uDelta)\le \epsilon q_{n-1,k+1}^{-1}$ and $\epsilon$ small enough, 
 the left-hand sides of both inequalities (21)  must vanish, since these are integers. 
 Then  the line $\Delta$, which contains the two points $ P_{n-1,k}$ and $P_{n-1,k+1}$,  coincides with  
 $ \Delta_{n,0}$, in contradiction with our assumptions. 
Therefore   the lower bounds
$$
L(\uDelta) >  \epsilon  q_{n-1,k+1}^{-1} = 16 \epsilon h_{n}^{-\sigma_{k+1}}
\and h  >  \epsilon \delta_{(n,k)} = {\epsilon\over 16} h_{n}^{\sigma_k -1} 
$$
hold,  so that 
$$
L(\uDelta) \gg \epsilon (h/\epsilon)^{ -\sigma_{k+1}/(\sigma_k -1)} \gg \epsilon (h/\epsilon)^{ -\mu},
$$
bounding $\sigma_{k+1} \le \sigma/\tau_0$ and $\sigma_k -1 \ge w-1$, since $k\ge 1$.

Consider now the case  $k=\ell'$. Then $h$ belongs to  the interval
$$
{\epsilon\over 16}  h_{n+1}^{\sigma - \tau_0} = {\epsilon q_{n,0}\over h_n}
< h  \le  {\epsilon  q_{n,1}\over h_{n+1}} = {\epsilon\over 16} h_{n+1}^{w-1}. \leqno {(22)}
$$
Arguing as in (21), we use here the single inequality
$$
 | r a_{n,0} + s b_{n,0} +  t c_{n,0} | 
 \ll  q_{n,0} L(\uDelta)  + \epsilon .
\leqno {(23)}
$$
If $\Delta$ does not pass  through the point $P_{n,0}$, the left-hand side of (23) is $\ge 1$,  and
noting that $\sigma \ge w\tau_0$, we obtain the required lower bound
$$
L(\uDelta) \gg q_{n,0}^{-1} \gg (h/\epsilon)^{-\sigma/(\sigma - \tau_0) } \gg (h/\epsilon)^{-\mu},
$$
provided $\epsilon$ is small enough.
 It remains to deal with  lines $\Delta$ containing the point  $P_{n,0}$. 
 Since $\Delta_{n+1,0}$ is the line joining $P_{n,0}$ and $P_{n+1,0}$,
 we may   apply   formula  (2)  to find 
$$
{ 1\over h}  \left | r { a_{n+1,0}\over c_{n+1,0}} + s { a_{n+1,0}\over c_{n+1,0}} +  t  \right |  = 
 d(P_{n+1,0}, \Delta) \gg\ll d(\Delta, \Delta_{n+1,0})d(P_{n,0},P_{n+1,0}).
\leqno{(24)}
$$
It readily follows from (22) and (1)  that 
$$
h>{\epsilon\over 16}  h_{n+1}^{\sigma - \tau_0}  \ge {\epsilon\over 16} h_{n+1}^{(w-1)\tau_0} 
\ge {\epsilon\over 16} h_{n+1}^{\tau_0} = {\epsilon\over 16} h_n.  
$$
  Accordingly,  we may define $k'$ 
as the largest integer $k \le \ell$ such that $ h > \epsilon h_{n,k} /16$. 
 Suppose first that $k'\le \ell -1$, so that $\epsilon h_{n,k'}/16 <  h \le \epsilon h_{n,k'+1}/16$.
 Then  Liouville inequality, applied to the  pencil of lines passing through  $P_{n,0}$, yields the lower bound
$$
d(\Delta ,\Delta_{n,k'}) \gg {q_{n,0} \over h_{n,k'} \, h} \gg \epsilon^{-1} { q_{n,0} \over h_{n,k'}h_{n ,k'+1}}.
$$
On the other hand, (12) gives the upper bound
$$
d(\Delta_{n,k'}, \Delta_{n+1,0}) \ll {q_{n,0} \over h_{n,k'} \, h_{n,k'+1}}.
$$
Using now the triangle inequality, the two above inequalities imply the lower bound
$$
d(\Delta , \Delta_{n+1,0}) \gg {q_{n,0} \over h_{n,k'} \, h}, \leqno {(25)}
$$
provided $\epsilon $ is small enough. Notice that (25) follows  directly from Liouville
inequality when $k'=\ell$. Next, combining (15), (24) and (25), we find the lower bound 
$$
 \left | r { a_{n+1,0}\over c_{n+1,0}} + s { b_{n+1,0}\over c_{n+1,0}} +  t  \right |
 \gg { h_{n+1} \over h_{n,k'} q_{n,1}}
 \gg h_{n+1}^{ -(w-1 + \tau_{k'})}. \leqno{(26)}
 $$
 On the other hand, it follows from (22) and (14) that
 $$
 h \, d(P_{n+1,0},\Theta) \ll \epsilon{q_{n,1}\over h_{n+1}} {h_{n+2}\over q_{n+1,0}q_{n+1,1}} 
 = 16 \epsilon h_{n+1}^{w-1-(w-1+\sigma)/\tau_0} \le 16 \epsilon h_{n+1}^{-w},
 $$
 noting that $\sigma \ge w\tau_0$ and $0< \tau_0 < 1$. 
  We can therefore replace in the left-hand side of (26) the coefficients $a_{n+1,0}/ c_{n+1,0}$
  and $b_{n+1,0}/ c_{n+1,0}$ by their limits $\alpha$ and $\beta$,  to  obtain the lower bound 
 $$
  | r \alpha + s \beta +  t  |  \gg h_{n+1}^{ -(w-1 + \tau_{k'})} \gg ( h/\epsilon)^{ -(w-1 + \tau_{k'})/\tau_{k'}}
  \gg (h/\epsilon)^{ -\mu} ,
  $$
 since $ h > \epsilon h_{n,k'}/16 =  \epsilon h_{n+1}^{\tau_{k'}}/16$. Fixing $\epsilon$ sufficiently small,
we have proved   the required lower bound $L(\uDelta)  \gg h^{-\mu}$  for any integer triple $\uDelta$
which is not a multiple of some $\uDelta_{n,k}$.
 
 \medskip
 
  Finally we prove the second part of Lemma 5.   Let $\uDelta$ be a non-zero  integer triple  with norm
   $h\le \epsilon h_{n+1} $.
The previous inequality  (23) remains valid. When $P_{n,0}$ does not lie on $\Delta$, we thus obtain  the stronger lower bound
$$
L(\uDelta) \gg q_{n,0}^{-1} \gg h_{n+1}^{-\sigma}  \gg h_{n+1}^{-w}.
$$
Suppose  now that $\Delta$  passes through $P_{n,0}$.  Notice that $\Delta$ cannot be equal to $\Delta_{n+1,0}$,
since the norm $h$ of $\uDelta$ is smaller than the height $H(\Delta_{n+1,0}) \ge  h_{n+1}$ of the line $\Delta_{n+1,0}$.  
Then,  we use Liouville inequality to bound from below
$$
d(\Delta , \Delta_{n+1,0}) \gg {q_{n,0}\over h \, h_{n+1}}.
$$
Taking again the argumentation leading to (26) with $k'=\ell$, we find  the lower bound
$$
 \left | r { a_{n+1,0}\over c_{n+1,0}} + s { b_{n+1,0}\over c_{n+1,0}} +  t  \right |
 \gg q_{n,1}^{-1} \gg h_{n+1}^{ -w}. \leqno{(27)}
 $$
As before, we may substitute in (27) the coordinates of the point $P_{n+1,0}$
by those of  $\Theta$, to obtain the required estimate
$$
  L(\uDelta) = | r \alpha + s \beta +  t  |  \gg h_{n+1}^{ -w} .
  $$
  \cqfd
 
\bigskip 

We   easily deduce from the assumptions (1)  that the strict upper bounds
$$
\displaylines{
\lambda = 
\max\Big( {1\over w-1 +\tau_0}, {\sigma -\tau_0\over \sigma }\Big) 
< {w-1\over \sigma}
\cr
\mu =  \max\Big( {\sigma\over (w-1 ) \tau_0}, {w-1+\tau_0\over \tau_0}\Big) = {w-1+\tau_0\over \tau_0}< {w-1+\tau_1\over \tau_0},
\cr }
$$
hold. Then Lemmas 4 and 5, together with (18),  show that the exponents of approximation
$\om(\Theta)$ and $\om({}^t\Theta)$ are reached respectively on  the set of integer triples
$(\uDelta_{n,k})_{n\ge 1, 0\le k \le \ell}$ and $(\uP_{n,k})_{n\ge 1, 0\le k \le \ell'}$. Now the estimates 
(16) and (17)  give the equalities
$$
\displaylines{
\om({}^t \Theta) =    \max_{0 \le k\le \ell' -1}\left({ \sigma_{k+1} -1\over \sigma_{k}}\right) = {w-1\over \sigma}
\cr
\om(\Theta)  =  \max_{0\le k \le \ell -1}\left({ w-1 + \tau_{k+1}\over \tau_k}\right) =  {w-1+ \tau_1\over \tau_0}.
\cr}
$$
The second parts of Lemmas 4 and 5 provide us with the upper bounds
$$
 \omc({}^t \Theta) \le {w-1\over w}
\and
\omc(\Theta) \le w.
$$
Taking into account the lower bounds (18), this concludes the proof of our proposition.
Notice  that  Lemma 5, together with (17),  yields a fine measure of linear independence over $\bQ$
of the numbers $1,\alpha,\beta$, which are obviously $\bQ$-linearly independent as required in the theorem.
 
 \section
{7. Infinite exponents}

The  basic construction  considered in Section 5 may  be greatly extended  by the introduction of variable exponents 
$\tau_{n,k}$ and $\sigma_{n,k}$,  depending  on $n$, instead of the fixed exponents $\tau_k$ and $\sigma_k$ occurring in (10). 
At each step $n$, we may also allow $\ell$ and $\ell'$ to vary (observe that Lemmas 1 and 2 are valid for any positive integers $\ell$ and $\ell'$). 
We take advantage of this flexibility
to complete the proof of the theorem in the remaining cases where $v= +\infty$. 
Our intention here is not to repeat the whole argumentation;  we briefly indicate below some 
specific choices of parameters $\tau_{n,k}$ and $\sigma_{n,k}$   leading   to any  quadruple of the form
$$
(+\infty, v', w,{w-1\over w}) \quad {\rm where} \quad  2\le w\le +\infty  , \,\, w-1 \le v' \le +\infty.
$$
Notice however that it might be
 useful to display more general constructions in order to compute 
the Hausdorff dimension of  subsets of points $\Theta\in \bR^2$ for which the quadruple of exponents
$\Omega(\Theta)$  belongs to various parts \footnote {(*)}{ As an example, the precise value of the Hausdorff dimension of the set 
$\{(\alpha,\beta)\in\bR^2; \omc(\alpha,\beta) \ge w\}$,  for a  given real 
number $w > 2$,  remains unknown. See [\Bak, \BuLaB , \Ry]  for estimates of that dimension in term of $w$. }
of $\bR^4$.

Let $w$ and $v'$ be real numbers with $w\ge 2$ and $v' \ge w-1$. Denote $\sigma =(w-1)/v'$ and
for any integer  $n > w/\sigma$, set
$$
\tau_{n,0} = {1\over n} , \quad \tau_{n,1}=1,\quad \sigma_{n,0} =\sigma , \quad \sigma_{n,1} = w .
$$
We first extend the increasing sequence $\sigma_{n,0} < \sigma_{n,1}$ using  an arithmetical progression with  $n$ terms
$$
w= \sigma_{n,1} < \cdots < \sigma_{n,n}= n\sigma,
$$
whose step $(n\sigma -w)/(n-1)$ is $<1$.  Then, the properties (7--9),  with $\ell=1$ and $\ell' =n$,
remain true  with  our present choice of parameters. The assumption $\sigma \le 1$ yields  the fundamental upper bound
$\sigma_{n,0}\le \tau_{n,0}+\tau_{n,1}$ occurring in (1). 
Next we
fix   two increasing sequences of positive real numbers  $(h_n)_{n >  w/\sigma }$ and $(q_{n,k})_{n>w/\sigma, 0 \le k \le n}$, 
satisfying  the recurrence relations 
$$
h_{n+1}= h_n^n =  h_n^{1/\tau_{n,0}}   , \qquad  q_{n,k} =  h_{n+1}^{\sigma_{n,k}}/16, \quad (0\le k\le n).
$$
The compatibility  relations $q_{n+1,0}= q_{n,n}$ hold for any $n$. Going  again through  the  construction described in Section 5, we obtain 
 a point $\Theta=(\alpha,\beta)$ with
$$
\eqalign{
\om(\Theta) =  &   \limsup_{n\rightarrow +\infty}\left({\sigma_{n,1} -1+ \tau_{n,1}\over \tau_{n,0}}\right)= +\infty,\quad
\cr
 \om({}^t\Theta) = &  \limsup_{n\rightarrow +\infty}\max_{0\le k\le n-1}\left( {\sigma_{n,k+1}-1\over \sigma_{n,k}}\right) 
  = {w-1\over \sigma} = {v'} ,
\cr
\omc(\Theta)= & \liminf_{n\rightarrow +\infty} \left({\sigma_{n,1} -1+ \tau_{n,1}\over \tau_{n,1}}\right)= w, 
\cr
\omc({}^t\Theta)= & \liminf_{n\rightarrow +\infty}\min_{0\le k\le n-1} \left({\sigma_{n,k+1}-1\over \sigma_{n,k+1}}\right)={w-1\over w}. 
\cr}
$$
We omit the details of the proof  which follows mutatis mutandis the same lines  as for the proposition.
Notice that Lemma 4 remains actually  valid with the exponent $\lambda =1$. 

When $v'=+\infty$, we make   use of sequences $(\sigma_{n,0})_n$ tending  to $0$. If $w\ge 2$ is a real number, take  
 $\ell= \ell' =1$ and set 
$$
\tau_{n,0} = {1\over n} , \quad \tau_{n,1}=1,\quad \sigma_{n,0} ={w \over n}, \quad \sigma_{n,1} = w
, \quad (n\ge  w). 
\leqno{(28)}
$$
Then,  we obtain a point $\Theta$ such that
$$
\om(\Theta) =\om({}^t\Theta)=+\infty
\and
\omc(\Theta)= w, \quad \omc({}^t\Theta) ={w-1\over w} .
$$
If moreover  $w=+\infty$,  substitute (for example)  $\sqrt{n}$   for  $w$   in the formulas (28). 
In that case,  the construction  produces   a point $\Theta=(\alpha,\beta)$ with $1,\alpha,\beta$ linearly independent over
$\bQ$, such that
$$
 \om(\Theta) =  \om({}^t\Theta)=\omc(\Theta)=+\infty \and  \omc({}^t\Theta) =1 .
$$

\vskip 2cm

\centerline{\bf References }

\vskip 5mm

\item{[\Bak]}
R. C. Baker,
{\it Singular $n$-tuples and Hausdorff dimension. II},
Math. Proc. Cambridge Phil. Soc. 111 (1992), 577--584.

\item{[\BoVa]}
E. Bombieri and J. Vaaler,
{\it On Siegel's Lemma},
Invent. Math. 73 (1983), 11--32.

\item{[\BuLaA]}
Y. Bugeaud and M. Laurent,
{\it Exponents of Diophantine Approximation and
Sturmian Continued Fractions}, Ann. Inst. Fourier, 55, 3 (2005), 
773--804.

\item{[\BuLaB]}
Y. Bugeaud and M. Laurent,
{\it On exponents of homogeneous and inhomogeneous Diophantine approximation}, 
to appear in the Moscow Mathematical Journal.

\item{[\BuLaC]}
Y. Bugeaud and M. Laurent,
{\it On exponents of Diophantine Approximation}, 
Preprint.

\item{[\DaSc]}
          {H. Davenport and W. M. Schmidt},
{\it Approximation to real numbers by
algebraic integers}, Acta Arith. {15} (1969), 393--416.

\item{[\Fis]}
          { S. Fischler},
{\it Spectres pour l'approximation d'un nombre r\'eel et de son carr\'e}, C. R. Acad. Sci. Paris, Ser. I  {339} (2004), 679--682.

\item{[\JarA]}
V. Jarn\'\i k,
{\it \"Uber ein Satz von A. Khintchine},
Pr\'ace Mat.-Fiz. 43 (1935), 1--16.

\item{[\JarB]}
V. Jarn\'\i k,
{\it O simult\'ann\i\' ch diofantick\'ych approximac\i\' ch},
Rozpravy T\'r. \v Cesk\'e Akad 45, c. 19 (1936), 16 p.

\item{[\JarC]}
V. Jarn\'\i k,
{\it \"Uber ein Satz von A. Khintchine, 2. Mitteilung},
Acta Arith.  2  (1936), 1--22.

\item{[\JarD]}
V. Jarn\'\i k,
{\it Zum Khintchineschen ``\"Ubertragungssatz''},
Trav. Inst. Math. Tbilissi 3 (1938), 193--212.

\item{[\JarE]}
V. Jarn\'\i k,
{\it Une remarque sur les approximations diophantiennes lin\'eaires},
Acta Sci. Math. Szeged 12 (1950), 82--86.

\item{[\JarF]}
V. Jarn\'\i k,
{\it Contribution \`a la th\'eorie des approximations diophantiennes
lin\'eaires et homog\`enes},
Czechoslovak Math. J. 4 (1954), 330--353 (in Russian, French summary).

\item{[\Khi]}
A. Ya. Khintchine,
{\it \"Uber eine Klasse linearer diophantischer Approximationen},
Rendiconti Circ. Mat. Palermo 50 (1926), 170--195.

\item{[\Lag]}
          J. C. Lagarias,
{\it Best Diophantine approximations to a set of linear forms},
J. Austral. Math. Soc. Ser. A  {34} (1983), 114--122.

\item{[\RoyA]}
          D. Roy,
{\it Approximation to real numbers by cubic algebraic numbers, I},
Proc. London Math. Soc. 88 (2004), 42--62.

\item{[\RoyB]}
          D. Roy,
{\it Diophantine approximation in small degree},
Number Theory, C.R.M. Proc. Lecture Notes 36 (2004), 269--285.

\item{[\RoyC]}
          D. Roy,
{\it On two exponents of approximation related to a real number and its square},
Canad. J. Math. To appear.

\item{[\Ry]}
B. P. Rynne,
{\it A lower bound for the Hausdorff dimension of
sets of singular $n$-tuples},
Math. Proc. Cambridge Phil. Soc. 107 (1990), 387--394.

\item{[\Sch]}
W. M. Schmidt,
{\it On heights of algebraic subspaces and Diophantine approximations},
Ann. Math.  85 (1967), 430--472.

\bigskip

\noindent   {Michel LAURENT}

\noindent 
{Institut de Math\'ematiques de Luminy}

\noindent 
{C.N.R.S. -  U.M.R. 6206 - case 907}

\noindent       {163, avenue de Luminy}

\noindent 
{13288 MARSEILLE CEDEX 9  (FRANCE)}

\noindent 
{\hbox{\tt laurent@iml.univ-mrs.fr}}

\bye